\documentclass[a4paper,11pt]{article}
\usepackage{graphicx}
\usepackage{amsfonts}
\usepackage{amsmath}
\usepackage[utf8]{inputenc}
\usepackage{url}
\usepackage{color}
\usepackage{mathtools}

\newcommand{\Rea}{\mbox{$\mathbb{R}$}}

\newcommand{\lip}{{\rm Lip}}

\newtheorem {Proposition}{Proposition} [section]
\newtheorem {Lemma}[Proposition] {Lemma}
\newtheorem {Theorem}[Proposition]{Theorem}
\newtheorem {Corollary}[Proposition]{Corollary}
\newtheorem {Remark}[Proposition]{Remark}
\newtheorem {Example}{Example}[section]

\begin{document}
\title{\sc Box-constrained monotone $L_\infty$-approximations to Lipschitz regularizations, with applications  to robust testing.
\footnote{Research partially supported by 
FEDER, Spanish Ministerio de Econom\'{\i}a y Competitividad, grant MTM2017-86061-C2-1-P and Junta de Castilla y Le\'on, grants VA005P17 and VA002G18.}}

\author{Eustasio del Barrio,    Hristo Inouzhe and Carlos Matr\'an\\
\textit{Departamento de Estad\'{\i}stica e Investigaci\'on Operativa and IMUVA,}\\
\textit{Universidad de Valladolid. SPAIN} }
\maketitle

\begin{abstract}
Tests of fit to exact models in statistical analysis often lead to rejections even when the model
is a useful approximate description of the random generator of the data. Among possible relaxations of a fixed model, 
the one defined by contamination neighbourhoods, namely, $\mathcal{V}_\alpha(P_0)=\{(1-\alpha)P_0+\alpha Q: Q \in \mathcal{P}\}$, where 
$\mathcal{P}$ is the set of all probabilities in the sample space, has received much attention, from its central role in Robust Statistics.
For probabilities on the real line, consistent tests of fit to $\mathcal{V}_\alpha(P_0)$ can be based on $d_K(P_0,R_\alpha(P))$,
the minimal Kolmogorov distance between $P_0$ and the set of trimmings of $P$, $R_\alpha(P)=\big\{\tilde P\in\mathcal{P}:\tilde P\ll P,\,{\textstyle \frac{d\tilde P}{dP}\leq\frac{1}{1-\alpha}}\, P\text{-a.s.}\big\}$.
We show that this functional admits equivalent formulations in terms of, either best approximation in uniform norm by $L$-Lipschitz functions satisfying a box constraint, or as the best monotone approximation in uniform norm to the $L$-Lipschitz regularization, which is seen to be expressable in terms of the average of the Pasch-Hausdorff envelopes. This representation for the solution of the variational problem allows to obtain
results showing stability of the functional $d_K(P_0,R_\alpha(P))$, as well as directional differentiability, providing the basis for a Central Limit Theorem for that functional. 
\end{abstract}


\noindent
{\it A.M.S. classification:} {\sc Primary:}  49J30. {\sc Secondary:} 26A16, 62G35, 41A29.

\section{Introduction.}
A repeated joker phrase in Statistics says that all models are wrong, but some  are useful. This celebrated aphorism,  attributed to the statistician G. Box,  on the one hand cautions that all models are approximations, while, on the other, stresses  the usefulness of good approximate models. Here, approximation should be interpreted, in words of Davies \cite{DaviesK}, as ``some formal admission of the fact that the statistical models are not true representations of the data". From this perspective,  within  the research objectives of Mathematical Statistics, it becomes natural the permanent  interest in the design and analysis of well-behaved procedures under small variations in the model. This includes the reconsideration of excessively restrictive concepts in Statistics, such as exact fit to models (say in  homogeneity, regression or time series settings). The interest is not exact equality, but only  ``similarity" or, alternatively, to find a ``relevant" difference. Also notice that this concept is  of great relevance in some applications, such as bioequivalence in Biostatistics (see, for example, \cite{Wellek}). Some recent references  sharing this spirit are (\cite{Munk}, \cite{modelAdec}, \cite{Dette}, \cite{Dette2}, \cite{Dette3}, \cite{Alvarez-Esteban2016}, \cite{Alvarez-Esteban2017}). 
That is also the perspective of our recent work \cite{Hristo}, while the present  paper addresses the mathematical bases giving  support to the approach.

 Let us begin  some historical notes on approximate model checking. A pioneer work in that sense is   \cite{hodges1954testing}. There, Hodges and Lehmann pointed out that ``when testing statistical hypotheses, we usually do not wish to take the action of rejection unless the hypothesis being tested is false to an extent sufficient to matter''. This fact led them to establish a distinction between statistical significance and material significance in hypotheses testing and to suggest modifications of the customary tests, in order to test for the absence of material significance. Their approach was based on assuming a distance in the parametric space and to allow some little deviation in the null hypothesis of the model. 
 
Ten years later, in his seminal paper \cite{huber1964robust}, Huber introduced the contamination neighbourhood of a probability,  namely, 
\begin{equation}
\label{def_cont_neigh}
\mathcal{V}_\alpha(P_0)=\{(1-\alpha)P_0+\alpha Q: Q\in\mathcal{P}\},
\end{equation}
where $\mathcal{P}$ is the set of all probability distributions in the space. Thus the probabilities in the neighbourhood are  mixtures of level $\alpha$ of $P_0$ with other probabilities. Although it can be defined in a wholly general setting,  throughout the paper $\mathcal{P}$ will be the set of probabilities on the (Borel) sets, $\beta$, of the real 
line $\mathbb{R}$). In this way, given an ``ideal" model $P_0$,  the vicinity includes  those probabilities which are distorted versions of the model through gross errors: given a particular value $\alpha_0\in [0,1)$, a probability $P$ in $\mathcal{V}_{\alpha_0}(P_0)$ would generate samples with an approximate   $(1-\alpha_0)\times100$  percentage of data coming from $P_0$.

Contamination neighbourhoods become one of the very basis of  Robust Statistics: a general attempt to provide methods with good performance when there are small departures from the assumed model.  Not surprisingly, its simple interpretation in terms of mixtures, motivated their use in different settings. In particular, Rudas et al \cite{rudas1994new} introduced a new index of fit in the framework of contingency tables.  Now the goal is to evaluate how well the contaminated version of the model describes the data, so statistical evidence of a ``small $\alpha$" should be considered as almost agreement with the model. The reconsideration of the problem in \cite{parMult}, also in the multinomial setting, allowed little deviations of the model that are measured by the Kullback--Leibler divergence. 

{The setup considered in both \cite{rudas1994new} or \cite{parMult} is constrained to the case when $P_0$ is a probability with a finite support. 
It should be noted at this point that testing fit to a neighbourhood of a fixed probability is not always a feasible task, depending on the metric or divergence which determines the neighbourhoods. Barron \cite{uecBa} considered the problem testing fit to approximate 
models and argued, while consistent tests were available for some weak metrics, it would
desirable that from the statistical assessment that $P$ and $P_0$ are close in a certain metric one could conclude that $P$ and $P_0$ are close from every point of view. In plain words, he advocated for the use of strong metrics, such as the total variation metric $d_{TV}(P,P_0)=\sup_{A\in\mathcal{A}}|P(A)-P_0|$, where $\mathcal{A}$ denotes the class of all measurable sets. Unfortunately, he also showed that if the probability $P_0$ is not discrete then there is no uniformly consistent test of fit to $P_0$ against alternatives at a certain distance in $d_{TV}$ and the same result remains true if the alternatives are bounded from $P_0$ in a distance or divergence that dominates the total variation metric. 
With these limitations in mind one may wonder if uniformly consistent testing to a meaningful relaxation of total variation neighbourhoods is possible beyond the discrete setting. In fact, contamination neighbourhoods are related to total variation neighbourhoods: $d_{TV}(P,Q)\leq \alpha$ if and only if there exists a probability 
$P_0$ such that $P\in\mathcal{V}_\alpha(P_0)$ and $Q\in\mathcal{V}_\alpha(P_0)$
(see \cite{simTrim}).

In \cite{Hristo} we showed that it is possible indeed to build a uniformly consistent test of fit to a contamination neighbourhood against increasingly closer alternatives.
We addressed the problem through the dual approach of trimmed probabilities,
an idea that goes back at least to  \cite{Gordaliza}. A probability $\tilde P\in \mathbb R$ is said to be a trimming of level $\alpha\in[0,1)$ of $P$ whenever there exists a down-weighting function $w$ such that 
$0\leq w \leq 1$ and $\tilde P(B)=\frac 1 {1-\alpha}\int_Bw(x)P(dx)$ for all the sets $B\in\beta$. Equivalently, it must be absolutely continuous w.r.t. $P$, with Radon-Nykodim derivative bounded by $\frac{1}{1-\alpha}$. The set of $\alpha$-trimmings of the probability distribution $P$ will be denoted by $R_\alpha(P)$:
\begin{equation}
\label{def_alpha_trimm}
R_\alpha(P)=\big\{\tilde P\in\mathcal{P}:\tilde P\ll P,\,{\textstyle \frac{d\tilde P}{dP}\leq\frac{1}{1-\alpha}}\, P\text{-a.s.}\big\}.
\end{equation}
The key link between (\ref{def_cont_neigh}) and (\ref{def_alpha_trimm}), obtained in   \cite{optTrim}, is given by
\begin{equation}\label{equivalencia}
P\in \mathcal{V}_\alpha(P_0) \Longleftrightarrow P_0\in R_\alpha(P).
\end{equation}
This duality has been exploited for analysis of similarity between samples in a fully nonparametric context (\cite{simTrim}), or for the consideration of a testable almost stochastic dominance model (\cite{Alvarez-Esteban2016}, 
\cite{Alvarez-Esteban2017}). 
There is a subtle, but important consequence of the duality (\ref{equivalencia}). In a realistic statistical setting we do not know either the value $\alpha$ or the ``contaminated" distribution $P$ but we only have an approximation $\hat P$ to $P$ (usually  $\hat P$ will be the empirical measure associated to a data set), and our goal is to search for statistical evidence, based on $\hat P$, for or against the hypothesis $P\in \mathcal{V}_\alpha(P_0)$. It turns out that sets of trimmings are often well behaved with respect to some of the most useful metrics in Statistics, while contamination neighbourhoods are not.
If $d$ is a metric on $\mathcal P$ and $R_\alpha(P)$ is closed for $d$ then both conditions in (\ref{equivalencia}) are equivalent to 
\begin{equation}\label{otraequivalencia}
d(P_0,R_\alpha(P))=0. 
\end{equation}
With a suitable choice of $d$ we could also ensure that $d(P_0,R_\alpha(\hat P))$ is a
consistent estimator of $d(P_0,R_\alpha(P))$. The success  of this strategy will strongly depend on the suitability of the metric for this task. Our choice here is the Kolmogorov distance, $d_K$, that for two probabilities $P,Q\in \mathcal P$ is defined by  the $L_\infty$-distance between their distribution functions $F_P$ and  $F_Q$. Davies \cite{DaviesK} claims that {\it the Kolmogorov distance induces the natural topology for statistics. Firstly, random variables are generated at the level of
distribution functions \dots Secondly all diagnostic checks and model
validation techniques operate at the level of distribution functions and not at the level of density functions \dots}
We show in this work that for $d_K$ the equivalent characterization (\ref{otraequivalencia}) holds. Also, it follows from \cite{ratesConv} that it is posible to use the empirical version of $d(P_0,R_\alpha(P))$ to build a consistent estimator of $\alpha_0$, the minimal contamination level such that 
$P\in\mathcal{V}_\alpha(P_0)$ (see subsection 4.2 in \cite{Hristo}).
However, suitability of $d_K$ in this setting depends also on the feasibility of the generated procedures. In fact, some difficulties related this metric are well known, both for its mathematical analysis (lack of Fr\'echet or Hadamard differentiability) and for its computational aspect (lack of convergent algorithms).

The motivation of this work is to provide sound mathematical support to our approach in \cite{Hristo} focusing in tools for  diagnostics, comparison and validation of an approximate statistical model. We will show (see Lemma \ref{ref1}) that the minimal Kolmogorov distance to a set of trimmings can be represented in terms of a variational problem, as follows. We set $\Gamma=F_0(F^{-1})$, $F_0$ and $F$ being the distribution functions of $P_0$ and $P$. Then, with great generality,  the following identity holds:
\begin{equation}\label{reduction}
d_K(P_0,R_\alpha(P))=\min \{ \| h-\Gamma \|, h\in \mathcal C_\alpha\},
\end{equation}
where
\begin{equation}
\mathcal C_\alpha:=\{ h: [0,1] \to  [0,1] \mbox{ nondecreasing, with } h(0)=0, h(1)=1, \mbox{ and }  \| h\|_{\rm Lip}\leq 1 /{1-\alpha}\}.\label{theset}
\end{equation}
  Here, as will be used throughout, for any real valued mapping $f: \aleph \to \mathbb{R}$ defined on a metric space $(\aleph,d)$, with $\|f\|$ and 
  $\|f\|_{\rm Lip}$ we will denote the $L_\infty$ and the Lipschitz norms:
  \[\|f\|=\sup_{x\in \aleph}|f(x)|, \ \ \ \ \|f\|_{\rm Lip}=\sup_{x,y\in \aleph}\frac{|f(x)-f(y)|}{d(x,y)}.
  \]

The representation in (\ref{reduction}) translates the problem of best trimmed approximation in Kolmogorov distance into finding a useful expression for a best $L_\infty$-approximation to a monotone function by monotone, Lipschitz-continuous functions satisfying the boundary conditions $h(0)=0, h(1)=1$. We will show (see Theorem \ref{prop_disting_h}) that the solution to this problem can be expressed in terms of Pasch-Hausdorff envelopes (see \cite{Rockafellar}). We will also relate this process with the alternative way of obtaining  Ubhaya's monotone $L_\infty$-best approximation (see \cite{Ubhaya1,Ubhaya2}) to the Lipschitz regularization of the objective function.

There are two main implications of our analysis of the variational problem in (\ref{reduction}) in statistical applications. First, it proves the validity of a simple, fast algorithm introduced in \cite{Hristo} for the computation of the empirical estimator
$d_K(P_0,R_\alpha(\hat{P}))$. Additionally, we use it to prove a result on directional differentiability of the $L_\infty$-distance to the regularized version (see Corollary \ref{derivability}). The relevance of this type of results on directional differentiability  has been pointed out in \cite{Shapiro}, and recently highlighted in relation with statistical applications in \cite{Carcamo}. {In fact, these results provide the mathematical foundation allowing a Central Limit Theorem (see Theorem 4.1 in \cite{Hristo}), thus incoming statistical applications of the proposed methodology.}
We should note that, under the false-model paradigm, this Central Limit Theorem yielded some tools for comparing models or for determining the usefulness of particular models following lines  related to \cite{modelAdec},\cite{davies1995data} or \cite{DaviesCRC}. In particular we should highlight the applications in the False-Discovery-Rate  (FDR) setting (as considered e.g. in \cite{Genovese} or \cite{Meinshausen}). In \cite{Hristo} (see Section 5 there) we discuss on the applicability of our approach to that setting.

The rest of this work is organized as follows. In Section \ref{topologia} we will present some alternative characterization of the set $R_\alpha(P)$ as well as its main topological properties in the $L_\infty$ setting. We include a key result on the stability of the constrained regularizations (see Proposition \ref{propiedades})
as well as the announced variational representation (Lemma \ref{ref1}) and the solution
of the variational problem (Theorem \ref{prop_disting_h}). The proof of this result will follow from that of Theorem 3.3 in Section 3, which discusses best $L_\infty$ approximation by Lipschitz functions with box constraints. The key link here is that Pasch-Hausdorff envelopes preserve monotonicity. Under continuity (Theorem \ref{fcont}) we provide a more convenient representation of the minimal distance between a nondecreasing function and its best Lipsichtz approximation. Section \ref{Ubhaya_section} considers the related problem of best $L_\infty$ approximation by monotone functions with box constraints, generalizing Ubhaya's results.
Finally, Section \ref{Ubhaya_section} contains also the announced results on directional differentiability (Theorem \ref{asymptotic1}
and Corollary \ref{derivability}).

\section{The set of trimmings in the $L_\infty$-topological setting}\label{topologia}
Since probabilities on $(\mathbb{R},\beta)$ are determined by their  distribution functions (d.f.'s in the sequel) and (\ref{def_cont_neigh}) and (\ref{def_alpha_trimm}) can be equivalently stated in terms of the corresponding distribution functions, we will use the same notation $R_\alpha(F)$ and $\mathcal{V}_\alpha(F_0)$, with the same meanings as before, but 
defined in terms of distribution functions. On the other hand, the Kolmogorov distance between probabilities is defined just through the $L_\infty$-distance between the corresponding d.f.'s, but we will often keep the notation $d_K$ for this distance. 

The set $R_\alpha(F)$ can be also characterized, as shown in 
\cite{TrimmedCompar} (see also Proposition 2.2 in \cite{optTrim} for a more general result), in terms of the set  of 
$\alpha$-trimmed versions of the uniform probability $U(0,1)$. Notice that this set is just   $\mathcal{C}_\alpha$,  as defined in (\ref{theset}). The parameterization, obtained through the composition of the functions $h$ and $F$: $F_h=h\circ F$ gives
\begin{equation}\label{calpha}
R_\alpha(F)=\{F_h: h \in \mathcal{C}_\alpha\}.
\end{equation}
We note that, as a consequence, the ``trimmed Kolmogorov distance" from $F$ to $F_0$ is
$$d_K(F_0,R_\alpha(F)):=\inf_{\tilde{F}\in R_\alpha(F)}\|\tilde{F} -F_0\|=\inf_{h\in \mathcal{C}_\alpha}\|h\circ F -F_0\|.$$

The set $R_\alpha(F)$ is convex and also well behaved w.r.t.  weak convergence of probabilities and 
widely employed probability metrics (see Section 2 in \cite{optTrim}). 
We show next that this also holds for $d_K$.
\begin{Proposition}
\label{dkcompact}
For $\alpha \in (0,1)$ and distribution functions $F$, $F_0, F_1, F_2,G_1$ and $G_2$, we have:
\begin{itemize}
\item[(a)] $R_\alpha(F)$ is compact w.r.t. $d_K$.
\item[(b)] $d_K(F_0,R_\alpha(F))=\min_{\tilde{F}\in R_\alpha(F)}\|\tilde{F} -F_0\|=\min_{h\in \mathcal{C}_\alpha}\|h\circ F -F_0\|$.
\item[(c)] $
|d_K(G_1,R_\alpha(F_1))-d_K(G_2,R_\alpha(F_2))| \leq \|G_1-G_2\|+ {\textstyle \frac 1 {1-\alpha}}\|F_1-F_2\|.
$
\end{itemize}
\end{Proposition}
\medskip
\noindent \textbf{Proof.}
By the Ascoli-Arzelà Theorem, $\mathcal{C}_\alpha$ is a compact subset of the space of continuous functions on $[0,1]$ endowed
with the uniform norm. Hence, from any sequence of elements in $R_\alpha(F)$, say $\{h_n\circ F\}$ (recall (\ref{calpha})), we can extract
a uniformly convergent subsequence $h_{n_j}\to h_0 \in \mathcal{C}_\alpha$. But then, obviously, $h_{n_j}\circ F\to h_0\circ F$ in $d_K$, which proves (a).
Since, on the other hand,
$$\big|\|h_1\circ F-F_0\|-\|h_2\circ F-F_0\|\big|\leq \|h_1\circ F-h_2\circ F\|\leq \|h_1-h_2\|,$$
we see that the map $h\mapsto \|h\circ F-F_0\|$ is continuous and, consequently, it attains its minimum in $R_\alpha(F)$, as claimed in (b).
Finally, to check (c) we note that
\begin{eqnarray}\label{in1}
\lefteqn{\big|d_K(G_1,R_\alpha(F_1))-d_K(G_1,R_\alpha(F_2))\big|\leq  \sup_{h\in \mathcal{C}_\alpha} \big| \|G_1-h\circ F_1\|-\|G_1-h\circ F_2\| \big|}\hspace*{6cm}\\
\nonumber 
&\leq &\sup_{h\in C_\alpha}  \|h\circ F_1-h\circ F_2\|\leq {\textstyle \frac 1 {1-\alpha}} \|F_1-F_2\|
\end{eqnarray}
and 
\begin{eqnarray}\label{in2}
\big|d_K(G_1,R_\alpha(F_2))-d_K(G_2,R_\alpha(F_2))\big|\leq  \sup_{h\in \mathcal{C}_\alpha} \big| \|G_1-h\circ F_2\|-\|G_2-h\circ F_2\| \big|\leq
\|G_1-G_2\|.
\end{eqnarray}
Now, (\ref{in1}) and (\ref{in2}) yield (c).
\quad $\Box$
\vspace{5mm}

Proposition \ref{dkcompact} guarantees the existence of optimal $L_\infty$-approximations to every distribution function $F_0$ by $\alpha$-trimmed versions of $F$: 
\begin{equation}\label{minimizers}
\mbox{There exists } \ \tilde F\in  R_\alpha(F) \ \mbox{ such that } \ \|F_0-\tilde F\|= d_K(F_0,R_\alpha(F)). 
\end{equation}
It also shows, through (\ref{equivalencia}), that for $\alpha \in [0,1)$
\begin{equation}\label{equivalencia2}
F \in \mathcal{V}_\alpha(F_0) \ \mbox{ if and only if } d_K(F_0,R_\alpha(F))=0.
\end{equation}
Moreover, by convexity of $R_\alpha(F)$, the set of optimally trimmed versions of $F$ associated to problem (\ref{minimizers}) is also convex.
However, guarantying uniqueness of the minimizer (as it holds w.r.t. $L_2$- Wasserstein metric by Corollary 2.10 in \cite{optTrim}) is not possible here.

An additional consequence of Proposition \ref{dkcompact} is the continuity of $d_K(F_0,R_\alpha(F))$ in $F_0$ and $F$. We quote this and some additional facts
in our next result.

\begin{Proposition}
\label{propiedades}
For $\alpha \in [0,1)$, if $\{F_n\}$ and $F$ are 
d.f.'s such that $d_K(F_n,F)\to 0,$ then:
\begin{itemize}
\item[a)] for every $\tilde{F}\in R_\alpha(F),$  there exist $\tilde{F}_n\in R_\alpha(F_n), n\in \mathbb N$ such that $d_K(\tilde{F}_n,\tilde{F})\to 0.$
\item[b)] if $\tilde{F}_n \in R_\alpha(F_n), n\geq 1$, then there exists some $d_K$-convergent subsequence $\{\tilde{F}_{n_k}\}$. If $\tilde{F}$ is the limit of 
such a subsequence, necessarily $\tilde{F}\in  R_\alpha(F)$.
\item[c)] if, additionally,  $\{G_n\}$ and $G$ are d.f.'s  such that $d_K(G_n,G)\to 0,$ then $d_K(G_n,$ $R_\alpha(F_m)) \to d_K(G,R_\alpha(F))$ as $n,m \to \infty.$
\end{itemize}
\end{Proposition}
\medskip
\noindent \textbf{Proof.}
To prove a), since $\tilde{F}=h\circ F$, with $h\in \mathcal{C}_\alpha$, it suffices to consider 
$\tilde{F}_n:=h \circ F_n \in R_\alpha(F_n)$ and recall that $h$ is Lipschitz. 
For b), we write $\tilde{F}_n=h_n\circ F_n$ and argue as in the proof of Proposition \ref{dkcompact} to get a $d_K$-convergent 
subsequence $h_{n_k}\to h\in \mathcal{C}_\alpha$ from which we easily get $d_K(h_{n_k}\circ F_{n_k},h\circ F)\to 0.$ Finally c) is a direct 
consequence of Proposition \ref{dkcompact} (c).
\quad $\Box$
\vspace{5mm}

By Polya's uniform convergence theorem, if $F$ and $G$ are continuous and $\{F_n\}, \{G_n\}$ are sequences of d.f.'s  which, respectively, weakly 
converge to $F,G$, then they also converge  in the $d_K$-sense, therefore $d_K(G_n, R_\alpha(F_m))\to d_K(G, R_\alpha(F))$ holds. Also, a 
direct application of the Glivenko-Cantelli theorem and item c) above guarantee the following  strong consistency result.
\begin{Proposition}
\label{consistency}
Let $\alpha \in [0,1)$ and  $\{F_n\}$ be the sequence of empirical d.f.'s based on a 
sequence $\{X_n\}$ of independent random variables with distribution function $F$. If $\{G_n\}$ is any  sequence of distribution 
functions $d_K$-approximating the d.f. $G$ (i.e. $d_K(G_n,G)\to 0$), then:
$$d_K(G_n,R_\alpha(F_m))\to d_K(G,R_\alpha(F)),\ \mbox{ as } n,m \to \infty, \ \mbox{ with probability one.}$$
\end{Proposition}
\vspace{5mm}

 Given a d.f. $F$, we write $F^{-1}$ for the associated quantile function (or left continuous inverse function), namely, 
$F^{-1}(t):=\inf\{x | \ t\leq F(x)\}$. We recall that if $U$ is a uniformly distributed 
$U(0,1)$ random variable,  $F^{-1}(U)$ has d.f. $F$. Similarly, if $X$ has a continuous d.f. $F$,  the composed function $F_0\circ F^{-1}$ is  
the quantile function associated to
the r.v. $Y=F_0(X)$ .  As we show next, under some regularity assumptions $d_K(F_0,R_\alpha(F))$ can be expressed in terms of the function
$F_0\circ F^{-1}$.    We will see later the usefulness of this fact both for the asymptotic analysis and the practical computation of  $d_K(F_0,R_\alpha(F_n))$ 
when $F_n$ is an empirical d.f. based on a data sample $x_1,\dots,x_n$. Recall that then $F_n(x):=\frac 1 n \sum_{i=1}^nI_{(-\infty,x]}(x_i)$.

\begin{Lemma} \label{ref1}
Let $\alpha \in [0,1)$. If $F, F_0$ are continuous d.f.'s  and $F$ is additionally strictly increasing then
$$d_K(F_0,R_\alpha(F))=\min_{h\in \mathcal{C}_\alpha}\|h-F_0\circ F^{-1}\| \
\mbox{ and } \
d_K(F_0, R_\alpha(F_n))=\min_{h\in \mathcal{C}_\alpha}\|h-F_0 \circ F_n^{-1}\|.$$
\end{Lemma}
\medskip
\noindent \textbf{Proof.} For the first identity observe that
\begin{eqnarray*} 
\|h\circ F-F_0\|&=\sup_{x\in\Rea}|h(F(x))-F_0(x)|=\sup_{F(x)\in[0,1]}|h(F(x))-F_0(F^{-1}(F(x)))|\\
&=\sup_{t\in[0,1]}|h(t)-F_0(F^{-1}(t))|=\|h-F_0(F^{-1})\|.
\end{eqnarray*}
On the other hand, if $x_{(i)}, i=1,\dots,n,$ denote the ordered sample associated to $x_1,\dots,x_n$ (the same set of values but ordered in nondecreasing sense) and 
$$t_0 = 0,\quad t_i=\frac i n,\quad h_i=h(F_n(x_{(i)}))=h(t_i),\quad\text{and}\quad F_{0,i}=F_0(x_{(i)}), \quad 1\leq i\leq n.$$
Taking into account that  $h(F_n)$ and $F_0(F_n^{-1})$ are piecewise constant while $F_0$ and $h$ are non decreasing and continuous, we obtain
$$\|h(F_n)-F_0\|=\max_{1\leq i\leq n}\max\Big(F_{0,i}-h_{i-1},h_i-F_{0,i}\Big)=\|h-F_0(F_n^{-1})\|,$$
and the other identity follows  from Proposition \ref{dkcompact}, part (b).
\quad $\Box$
\vspace{5mm}

Our final result in this section provides a simple representation of $\min_{h\in \mathcal{C}_\alpha}\|h-F_0\circ F^{-1}\|$ (hence, of $d_K(F_0,R_\alpha(F))$). In this statement
{we assume that $\Gamma$ is a nondecreasing function taking values in $[0,1]$ (which is always the case if $\Gamma=F_0\circ F^{-1}$). Note that taking  right and left limits
at 0 and 1, respectively, we can assume that $F_0\circ F^{-1}$ is a nondecreasing (and left continuous) function from $[0,1]$ to $[0,1]$.

\begin{Theorem}
\label{prop_disting_h}
Let $\alpha \in [0,1)$. Assume $\Gamma:[0,1]\to [0,1]$ is a  nondecreasing function. Define $G(t)=\Gamma (t)-\frac t {1-\alpha}$, $U(t)= \sup_{t\leq s\leq 1} G(s)$, $L(t)=\inf_{0\leq s\leq t}G(s)$ and
$$\tilde{h}_\alpha(t)=\max\left(\min\left( {\textstyle \frac{U(t)+L(t)}{2} }, 0\right),{\textstyle  \frac{-\alpha}{1-\alpha}}\right).$$
Then, 
$$\min_{h\in \mathcal{C}_\alpha}\|h-\Gamma\|=\|\tilde{h}_\alpha-G\|.$$
\end{Theorem}
}
The proof of this result will be developed in Section \ref{Approximation}. In fact Theorem \ref{final} is just a rephrasing of this result. A look at that Theorem shows that $h_\alpha=\tilde{h}_\alpha+\frac{\cdot}{1-\alpha}$ is an element of $\mathcal{C}_\alpha$ such that
$\|h_\alpha-\Gamma\|=\min_{h\in \mathcal{C}_\alpha}\|h-\Gamma\|$, that is, $h_\alpha$ is an optimal trimming function in the sense described above. We recall that we do not claim uniqueness
of this minimizer, but this particular choice allows to compute $d_K(F_0,R_\alpha(F_n))$ for sample d.f.'s. Moreover,  Theorem \ref{prop_disting_h} even provides a simple way for the computation
of $d_K(F_0,R_\alpha(F))$ for theoretical distributions. Let us see an illustration of this use.

\begin{Example}[Trimmed Kolmogorov distances in the Gaussian model.]\label{EjemploNormal1}  {
Consider the case $F_0=\Phi$, $F=\Phi((\cdot-\mu)/\sigma)$, where $\Phi$ denotes the standard normal d.f., $\mu\in\mathbb{R}$ and $\sigma>0$.
Here we have $H^{-1}(t):=F_0\circ F^{-1}(t)=\Phi(\mu + \sigma\Phi^{-1}(t))$. 
We note that $w(t):=(H^{-1})'(t)\leq 1/(1-\alpha)$ if and only if $p(\Phi^{-1}(t))\geq 0$, where
\begin{equation}\label{px}
p(x)=(\sigma^2-1)x^2+2\mu\sigma x+ \mu^2 - 2\log((1-\alpha)\sigma).
\end{equation}
To avoid cumbersome computations we focus on the cases $\sigma=1$, $\mu\ne 0$ and $\mu=0$, $\sigma\ne 1$.

\smallskip
If $\sigma=1$ and $\mu>0$ then $p$ is linear with positive slope and we see that $w(t)\leq 1/(1-\alpha)$ if and only if
$t\geq t_0=\Phi\big(-\frac \mu 2+\frac 1 \mu \log(1-\alpha)\big)$. This means that $G(s)=H^{-1}(s)-s/(1-\alpha)$ is increasing in 
$[0,t_0]$ and decreasing in $[t_0,1]$. Since, $H^{-1}(0)=G(0)=0$, we have that, $\tilde{h}_\alpha (t) = 0$ for $t\in[0,t_1]$, 
where $t_1\in (t_0,1)$ is (the unique) solution to  $G(t_1)=0$, and $\tilde{h}_\alpha (t) = G(t)$ for $t\in [t_1,1]$. 
We conclude that $d_K(R_\alpha(N(\mu, 1)), N(0,1))=G(t_0)$. The case $\mu<0$ can be handled similarly to obtain
\begin{equation}\label{normaldominated}
d_K(R_\alpha(N(\mu, 1)), N(0,1))=
{\textstyle \Phi\big(\frac {|\mu|} 2+\frac 1 {|\mu|} \log(1-\alpha)\big)-
\frac 1 {1-\alpha}\Phi\big(-\frac {|\mu|} 2+\frac 1 {|\mu|} \log(1-\alpha)\big)},\quad \mu\ne 0.
\end{equation}

\smallskip
We focus now on the case $\mu=0$. If $\sigma^2< 1$, $p$ is a parabola with negative leading coefficient and discriminant $\Delta^2=8(\sigma^2-1)\log(\sigma(1-\alpha))>0$.
Hence, $p(x)$ is positive for $x\in (x_a,x_b)$ with $x_a= -\frac \Delta {2(1-\sigma^2)}$, $x_b=\frac \Delta {2(1-\sigma^2)}$.
Equivalently, $w(t)\leq 1/(1-\alpha)$ if and only if $t_a:=\Phi(x_a)\leq t\leq t_b:=\Phi(x_b)$. This means that $G$ is increasing 
in $[0,t_a)$, decreasing in $[t_a,t_b]$, increasing in $(t_b,1]$, $G(0)=0$ and $G(1)=-\alpha/(1-\alpha)$. Arguing as above,
we have $\tilde{h}_\alpha(t)=\min(G(t),0)$ for $0\leq t\leq \frac 1 2$, $\tilde{h}_\alpha(t)=\max(G(t),-\frac \alpha{1-\alpha})$ 
for $\frac 1 2\leq t\leq 1$, $\tilde{h}_\alpha(t_a)=0$ and $\tilde{h}_\alpha(t_b)=\frac{-\alpha}{1-\alpha}$. 
We conclude that $d_K(R_\alpha(N(\mu, \sigma^2)), N(0,1))=G(t_a)-\tilde{h}_\alpha(t_a)=\tilde{h}_\alpha(t_b)-G(t_b)$. 
Hence,
$$d_K(R_\alpha(N(0, \sigma^2)), N(0,1))={\textstyle\Phi\Big(\frac {-\sigma \frac \Delta 2 }{1-\sigma^2} \Big)-\frac 1 {1-\alpha}
\Phi\Big(\frac {- \frac \Delta 2 }{1-\sigma^2} \Big)},\quad \mbox{if }\sigma<1.
$$

\smallskip
If $1\leq \sigma \leq 1/(1-\alpha)$ then we have that $w(t)\leq 1/(1-\alpha)$ for all $t$ and $h_0=H^{-1}\in \mathcal{C}_\alpha$. In particular, $d_K(R_\alpha(N(0, \sigma^2)), N(0,1))=0$.

\smallskip
Finally, we consider the case $\sigma > 1/(1-\alpha)$. In this case
$p$ is positive for $x\notin [x_a,x_b]$ with $x_a= -\frac \Delta {2(\sigma^2-1)}$, $x_b=\frac \Delta {2(\sigma^2-1)}$. 
This means that $(H^{-1})'(t)>\frac 1 {1-\alpha}$ for $t\in (t_a,t_b)$ with $t_a=\Phi(x_a),
t_b=\Phi(x_b)$.
Therefore, $G$ is decreasing in $[0,t_a)$, increasing in $[t_a,t_b]$, decreasing in $(t_b,1]$, $G(0)=0$ and $G(1)=-\alpha/(1-\alpha)$.
Hence, $\tilde{h}_\alpha(t)=\max(G(t),\frac{G(t)+G(t_b)}2)$, $0\leq t\leq t_a$, 
$\tilde{h}_\alpha(t)=\frac{G(t_a)+G(t_b)}2$, $t_a\leq t\leq t_b$, $\tilde{h}_\alpha(t)=\min(G(t),\frac{G(t_a)+G(t)}2)$, $t_b\leq t\leq 1$. In particular,  
$d_K(R_\alpha(N(0, \sigma^2)), N(0,1))=\tilde{h}_\alpha(t_a)-G(t_a)=G(t_b)-\tilde{h}_\alpha(t_b)=\frac 1 2(G(t_b)-G(t_a))$, that is,  
$$d_K(R_\alpha(N(0, \sigma^2)), N(0,1))={\textstyle\Phi\Big(\frac {\sigma \frac \Delta 2 }{\sigma^2-1} \Big)-\frac {\Phi\Big(\frac { \frac \Delta 2 }{\sigma^2-1} \Big)-\frac \alpha 2}{1-\alpha}
},\quad \mbox{if }\sigma>{\textstyle \frac 1 {1-\alpha}}.
$$

\hfill $\Box$
}
\end{Example}

\section{Best $L_\infty$-approximations by Lipschitz-continuous functions with box constraints}\label{Approximation}
In this section we refresh the notation. The role of $1/(1-\alpha)$ will be played now by a generic Lipschitz constant $L$; our $\Gamma$ will be substituted by a  bounded function $f:\aleph \to \mathbb R$, where $(\aleph,d)$ is (at least at the beginning) a general metric space, while we maintain $[0,1]$ as the range of values. We will also use the notation $x\lor y$ (resp. $x\land y$) for the maximum (resp. minimum) of both numbers (or functions). Regarding the Lipschitz norm, recall the trivial inequalities
\begin{equation}\label{desigL}
\|f\land g\|_\lip, \|f\lor g\|_\lip \leq \|f\|_\lip \lor \|g\|_\lip.
\end{equation}

The first lemma collects some basic properties on the role of the Pasch-Hausdorff envelopes of a function to obtain a  Lipschitz-continuous best $L_\infty$-approximation with constrained Lipschitz constant. For the sake of completeness, we will also include a simple proof.

\begin{Lemma}\label{primero}
For a function $f: \aleph\to [0,1]$, given a constant $L\geq 0,$ let us consider
$$f_{L,1}(x):=\inf_{y\in \aleph}(f(y)+Ld(x,y)), \ \ \ f_{L,2}(x):=\sup_{y\in \aleph}(f(y)-Ld(x,y)).$$
\begin{itemize}
\item[(i)] This defines functions $f_{L,1}, f_{L,2}: \aleph \to \mathbb R$ such that
$0\leq f_{L,1}\leq f_{L,2}\leq 1.$
\item[(ii)] $f_{L,1}$ is the pointwise largest function $g: \aleph \to \mathbb R$ satisfying $g\leq f$ and $\|g\|_\lip\leq L$. Likewise $f_{L,2}$ is the pointwise smallest function $g: \aleph \to \mathbb R$ satisfying $g\geq f$ and $\|g\|_\lip\leq L$.
\item[(iii)] The average $f_L:=(f_{L,1}+ f_{L,2})/2$ satisfies $\|f_L\|_\lip \leq L$ and
$$\|g-f\|\geq\|f_L-f\|=\|f_{L,2}-f_{L,1}\|$$
for any function $g:\aleph\to\mathbb R$ such that $\|g\|_\lip \leq L.$
\end{itemize}
\end{Lemma}
\medskip
\noindent \textbf{Proof.}
Part (i) follows directly from the definitions of $f_{L,1}$ and $f_{L,2}$, because, for every $x\in \aleph$:
$$\inf_{y\in \aleph}f(y)\leq f_{L,1}(x) \leq f(x)+Ld(x,x)=f(x)= f(x)-Ld(x,x)\leq f_{L,2}(x)\leq\sup_{y\in \aleph}f(y).$$

To address part (ii) observe that, for arbitrary $x_1,x_2,y\in \aleph$, the triangle inequality for the distance implies 
$|Ld(x_1,y)-Ld(x_2,y)|\leq Ld(x_1,x_2),$  leading to the inequalities
$$|f_{L,j}(x_2)-f_{L,j}(x_1)|\leq Ld(x_1,x_2) \ \ \mbox{ for } j=1,2,$$
thus to $\|f_{L,j}\|_\lip \leq L, j=1,2.$ Now, if $g:\aleph \to \mathbb R$ satisfies $g\leq f$ and $\|g\|_\lip \leq L,$ then for  $x,y\in \aleph$:
$g(x)\leq g(y)+Ld(x,y)$ with equality if $x=y$. Hence
$$g(x)=\inf_{y\in \aleph}(g(y)+Ld(x,y))\leq \inf_{y\in \aleph}(f(y)+Ld(x,y))=f_{L,1}(x).$$
Analogously, it follows from $g\geq f$ and $\|g\|_\lip \leq L$ that $g\geq f_{L,2},$ proving (ii). 

As to part (iii), let $\epsilon:=\|g-f\|.$ Then $\|g\pm \epsilon\|_\lip =\|g\|_\lip$ and $g-\epsilon\leq f\leq g+\epsilon.$ Consequently, by part (ii),
$$g-\epsilon\leq f_{L,1}\leq f\leq f_{L,2}\leq g+\epsilon$$
This implies that
$$|f_L-f|=(f-f_L)\lor(f_L-f)\leq (f_{L,2}-f_L)\lor(f_L-f_{L,1})=\frac {f_{L,2}-f_{L,1}}2\leq\epsilon,$$
whence
$$\|f_L-f\|\leq\frac{\|f_{L,2}-f_{L,1}\|}2\leq\|g-f\|.$$
Since $\|f_L\|_\lip \leq \|f_{L,1}\|_\lip/2+\|f_{L,2}\|_\lip/2\leq L,$ taking $g=f_L$ gives the announced equality $\|f_L-f\|=\|f_{L,2}-f_{L,1}\|/2.$
\quad $\Box$
\vspace{5mm}

When $\aleph$ is a real interval and $f$ is non-decreasing, the functions $f_{L,1}$ and $f_{L,2}$ in Lemma \ref{primero} share also that property and can be alternatively expressed in terms of the the Ubhaya's monotone envelopes of the function $f(x)-Lx$. This is the content of the following lemma.

\begin{Lemma}\label{segundo}
Let $\aleph$ be a real interval, equipped with the usual distance $d(x,y)=|x-y|.$ If $f:\aleph \to [0,1]$ is non-decreasing, then the functions $f_{L,1}, f_{L,2}$ in Lemma \ref{primero} are non-decreasing too, and for arbitrary $x\in \aleph$ and $j=1,2,$
$$f_{L,j}(x)=\gamma_{L,j}(x)+Lx,$$
where $\gamma_{L,j}, j=1,2$ are the non-increasing functions
$$\gamma_{L,1}(x):=\inf_{y\in \aleph: y\leq x}(f(y)-Ly) \ \ \mbox{ and } \ \  \gamma_{L,2}(x):=\sup_{y\in \aleph: y\geq x}(f(y)-Ly).$$
In particular,
\begin{equation}\label{varios}
\|f_{L,2}-f_{L,1}\| = \|\gamma_{L,2}-\gamma_{L,1}\|=\sup_{y,x\in \aleph: y\leq x}(f(x)-f(y)-L(x-y)).
\end{equation}
\end{Lemma}
\medskip
\noindent \textbf{Proof.}
The representations of $f_{L,1}$ and $f_{L,2}$ in terms of $\gamma_{L,1}$ and $\gamma_{L,2}$ follow from the fact that for arbitrary $x,y\in \aleph,$
\begin{equation*}
   f(y)+Ld(x,y)
    \begin{cases*}
      =f(y)+L(x-y)=f(y)-Ly+Lx & if $y\leq x$ \\
      \geq f(x)=f(x)-Lx+Lx        & if $y\geq x,$
    \end{cases*}
  \end{equation*}
\begin{equation*}
   f(y)-Ld(x,y)
    \begin{cases*}
      =f(y)-L(y-x)=f(y)-Ly+Lx & if $y\geq x$ \\
      \leq f(x)=f(x)-Lx+Lx        & if $y\leq x,$
    \end{cases*}
  \end{equation*}
where the inequalities follow from  $f$ being non-decreasing. Note that both functions $\gamma_{L,1}$ and $\gamma_{L,2}$ are non-increasing, but adding the term $Lx$ to them leads to non-decreasing functions: For $x_1,x_2\in \aleph$ with $x_1<x_2,$ isotonicity of $f$ implies that
\begin{eqnarray*}
f_{L,2}(x_1) &=& \sup_{y\geq x_2}(f(y)-Ly+Lx_1)\lor \sup_{x_1\leq y \leq x_2}(f(y)-Ly+Lx_1)\\
 &\leq& (f_{L,2}(x_2)-Lx_2+Lx_1)\lor f(x_2)\\
 &\leq& f_{L,2}(x_2),
\end{eqnarray*}
and
\begin{eqnarray*}
f_{L,1}(x_2) &=& \inf_{y\leq x_1}(f(y)-Ly+Lx_2)\land \sup_{x_1\leq y \leq x_2}(f(y)-Ly+Lx_2)\\
&\geq& (f_{L,1}(x_2)+Lx_2-Lx_1)\land f(x_1)\\
&\geq& f_{L,1}(x_1),
\end{eqnarray*}
because $f_{L,1}\leq f \leq f_{L,2}.$
\quad $\Box$
\vspace{5mm}

Finally, let us include in the problem the boundary restrictions.

\begin{Theorem}\label{final}
Let $f: [0,1]\to [0,1]$ be non-decreasing. For $L\geq 1$ consider the function
\begin{eqnarray*}
\tilde f_L(x) &:=&(f_L(x)\lor(1-L+Lx))\land Lx \\
&=&((\gamma_L(x)\lor (1-L))\land 0)+Lx,
\end{eqnarray*}
where $\gamma_L:=(\gamma_{L,1}+\gamma_{L,2})/2,$ and $f_L, \gamma_{L,1}, \gamma_{L,2}$ are defined as in Lemmas \ref{primero} and \ref{segundo}. Then $\tilde f_L: [0,1] \to \mathbb R$ is non-decreasing and verifies $\tilde f_L(0)=0$ and $\tilde f_L(1)=1$ and $\|\tilde f_L\|_\lip \leq L,$ and for arbitrary functions $g: [0,1] \to \mathbb R$ with $g(0)=0$ and $g(1)=1$ and $\|g\|_\lip \leq L,$
\begin{eqnarray}\nonumber
\|g-f\| &\geq& \|\tilde f_L-f\| \\
&=& \max \Big\{ f_{L,2}(0),1-f_{L,1}(1), \sup_{0\leq y \leq x \leq 1} (f(x)-f(y)-L(x-y))/2\Big\} \label{minvalue}
\end{eqnarray}
\end{Theorem}
\medskip
\noindent \textbf{Proof.} 
  Let us begin noting that both expressions for $\tilde f_L$ are trivially equivalent from the relations between $\gamma_{L,j}$ and $f_{L,j}.$

That $\tilde f_L$ verifies the required properties easily follows from the preceding lemmas (recall also inequalities (\ref{desigL})).   
Let then $g: [0,1] \to \mathbb R$ with $\|g\|_\lip \leq L$. Also by the precedent lemmas,
$$\|g-f\| \geq \|f_L-f\|= \sup_{0\leq y \leq x \leq 1}(f(x)-f(y)-L(x-y))/2.$$
Under the additional constraint that $g(0)=0$, for arbitrary $x\in [0,1],$
$$f(x)-g(x)=f(x)-(g(x)-g(0))\geq f(x)-Lx,$$
whence
$$\|g-f\|\geq \sup_{0\leq x \leq 1}(f(x)-Lx)=f_{L,2}(0).$$
Analogously, the additional constraint $g(1)=1$ implies that
$$f(x)-g(x)=f(x)+(g(1)-g(x))-1 \leq f(x)+ L(1-x)-1,$$
whence
$$-\|g-f\| \leq \inf_{0\leq x \leq 1}(f(x)+L(1-x))-1=f_{L,1}(1)-1.$$
These considerations show that for any function $g: [0,1] \to \mathbb R$ verifying the conditions $g(0)=0,$ $g(1)=1$ and $\|g\|_\lip \leq L,$
$$\|g-f\| \geq \|f_L-f\| \lor f_{L,2}(0)\lor (1-f_{L,1}(1)).$$

The function $\tilde f_L$ satisfies the previous constraints on $g$, too, so
$$\| \tilde f_L-f\| \geq \|f_L-f \|\lor f_{L,2}(0)\lor (1-f_{L,1}(1)).$$
It remains to prove the reverse inequality. For $x\in [0,1]$, we have to distinguish three cases: If $1-L+Lx\leq f_L(x)\leq Lx,$ then $\tilde f_L(x)=f_L(x),$ so $|\tilde f_L(x)-f(x)|\leq \|f_L-f\|.$ If $f_L(x)>Lx,$ then $\tilde f_L(x)=Lx,$ and
\begin{equation*}
   f(x)-\tilde f_L(x)
    \begin{cases*}
      =f(x)-Lx\leq f_{L,2}(0), \\
      > f(x)-f_L(x)\geq -\|f_L-f\|.
    \end{cases*}
  \end{equation*}
  Similarly, if $f_L(x)<1-L+Lx,$ then $\tilde f_L(x)=Lx,$ and
  \begin{equation*}
   f(x)-\tilde f_L(x)
    \begin{cases*}
      =f(x)+L(1-x)-1\geq f_{L,1}(1)-1, \\
      < f(x)-f_L(x)\leq \|f_L-f\|.
    \end{cases*}
  \end{equation*}
\quad $\Box$
\vspace{5mm}

In the case, considered in Theorem \ref{final}, of a non-decreasing function $f$, since the functions $f_{L,j}$ are absolutely continuous and the relations $\gamma_{L,j}=f_{L,j}-Lx$ hold, all the functions $f_L,  \gamma_L, \gamma_{L,j}$ are absolutely continuous so  $\{\gamma_L\leq1-L\}, \{\gamma_L\geq0\}, \{\gamma_L\in [1-L,0]\}$ are compact sets and continuous functions attain their maximum values on these sets. This allows to  get  alternative expressions for (\ref{minvalue}) as given in the following theorem. We note that here and throughout  we use the convention that the $\max$ over an empty set equals $-\infty$.

\begin{Theorem}\label{fcont} Let $f: [0,1]\to [0,1]$ be non-decreasing and continuous and assume the notation  in Theorem \ref{final}. Then the following  alternative expressions for (\ref{minvalue}) hold:
\begin{eqnarray}\label{alter}
\| f-\tilde f_L \|&=& \max \left( \max_{x\in \mathcal{T}_1}\left( f(x)-Lx\right),\max_{x\in \mathcal{T}_2}\left(1-L+Lx-f(x)\right),\frac 1 2 \max_{1-L\leq \gamma_L(x)\leq 0}\left(\gamma_{L,2}(x)- \gamma_{L,1}(x)\right)\right) \ \ \ \ \ \ \ \\ \label{alter2}
&=&\max \left( \max_{x\in \mathcal{T}_1}\left( f(x)-Lx\right),\max_{x\in \mathcal{T}_2}\left(1-L+Lx-f(x)\right),\frac 1 2 \max_{(y,x)\in\mathcal{T}_3}\left(f(x)- f(y)-L(x-y)\right)\right).  
\end{eqnarray}
Here, we used the notation $\mathcal{T}_1=\{x\in [0,1]  : \gamma_L(x)\geq 0\},$ $\mathcal{T}_2=\{x\in [0,1]  : \gamma_L(x)\leq 1-L\},$ $\mathcal{T}_3= \{(y,x): 0\leq y \leq x \leq 1, 1-L\leq \frac 1 2 (f(y)+f(x)-L(y+x))\leq 0\}.$
\end{Theorem}

Once we know Theorem \ref{final}, a proof of this result would take advantage of the fact that the right-hand side in (\ref{alter}) is upper bounded by the same expression with the unrestricted maxima, which, by (\ref{varios}) is just the right-hand side in (\ref{minvalue}) when $f$ is continuous. However, with some additional effort we can obtain a more general result that does not requires the monotonicity assumption on the objective function and opens a way to address the directional differentiability of the functional $f\to \|f-\tilde f_L\|$. Both goals will be carried through the following section.

\section{Best $L_\infty$-approximations by monotone functions with box constraints}\label{Ubhaya_section}

The following theorem gives appropriate characterizations of the best approximation of a  bounded function (in uniform norm) by monotone functions
with a box constraint. Without this constraint, best approximation by monotone functions in the $L_\infty$-norm has been considered
in \cite{Ubhaya1,Ubhaya2}, with results that cover the case $A=-\infty$, $B=\infty$ in Theorem \ref{representation1} below. Notice that this theorem, based on Ubhaya's envelopes, would also provide an (arguably more involved) alternative proof for  Theorem \ref{final}. Notice that the function $G$ plays the role of the transformed function, $f(x)-Lx$ (the difference of two nondecreasing functions) in the previous section, while the scope here is general.

\begin{Theorem}\label{representation1}
Assume $G:[0,1]\to \mathbb{R}$ is a bounded function and $-\infty\leq A\leq B\leq \infty$. Define
$U(x)=\sup_{x\leq y \leq 1} G(y)$, $L(x)=\inf_{0\leq y \leq x} G(y)$, $\bar{G}(x)=(L(x)+U(x))/2$ and 
$$\bar{G}_{A,B}(x)=\max(\min(\bar{G}(x),B),A).$$ Then $U, L, \bar{G}$ and $\bar{G}_{A,B}$ are nonincreasing,
$L(x)\leq G(x)\leq U(x)$
and for every nonincreasing $h:[0,1]\to [A,B]$ we have
\begin{equation}\label{bestappr}
\|G-\bar{G}_{A,B}\|\leq \|G-h\|.
\end{equation}
Furthermore, if $G$ is continuous then $U, L, \bar{G}$ and $\bar{G}_{A,B}$ are also continuous and
\begin{eqnarray}\nonumber
\|G-\bar{G}_{A,B}\|&=&\max\left(\max_{ \bar{G}(x)\geq B} (G(x)-B),\max_{\bar{G}(x)\leq A} (A-G(x)),{\textstyle \frac 1 2}\max_{A\leq \bar{G}(x)\leq B} (U(x)-L(x))\right)\\
\label{rep1}
&=&\max\left(\max_{x\in \mathcal{T}_1} (G(x)-B),\max_{x\in \mathcal{T}_2} (A-G(x)),{\textstyle \frac 1 2}\max_{(y,x)\in \mathcal{T}_3} (G(x)-G(y))\right),
\end{eqnarray}
where $\mathcal{T}_1=\{x\in [0,1]:\, \bar{G}(x)\geq B \}$, $\mathcal{T}_2=\{x\in [0,1]:\, \bar{G}(x)\leq A \}$ and $\mathcal{T}_3=\{(y,x):\, 0\leq y\leq x\leq 1, A\leq \frac 1 2 (G(y)+G(x))\leq B \}$.
\end{Theorem}

\medskip
\noindent \textbf{Proof.} The bounds $L(x)\leq G(x)\leq U(x)$ are obvious, and also the fact that $U$ and $L$ are nonincreasing (hence, also $\bar{G}$ and $\bar{G}_{A,B}$). 

\noindent
$\bullet$\hspace{3mm}Next, consider some nonincreasing  $h: [0,1]\to [A,B]$ and $x\in[0,1]$. Since $L(x)\leq G(x)\leq U(x)$, we have that $G(x)=\bar{G}(x)$ whenever $U(x)=L(x)$. Hence, if 
$U(x)=L(x)\in [A,B]$ we have $\bar{G}_{A,B}(x)=G(x)$
and, consequently, $$0=|\bar{G}_{A,B}(x)-G(x)|\leq \|h-G\|.$$

\noindent
$\bullet$\hspace{3mm}Obviously, $\bar{G}_{A,B}(x)=B$ if $U(x)=L(x)>B$ and we still have that $$|\bar{G}_{A,B}(x)-G(x)|\leq |h(x)-G(x)|\leq\|h-G\|$$ and similarly for 
the case $U(x)=L(x)<A$. 

\noindent
$\bullet$\hspace{3mm}It remains to deal with the case $U(x)>L(x)$. For every $\varepsilon>0$ there exist $x_a \in[0,x]$, $x_b\in[x,1]$ such that
$G(x_a)<L(x)+\varepsilon$ and $G(x_b)>U(x)-\varepsilon$. 
If $\bar{G}(x)>B$ then $\bar{G}_{A,B}(x)=B$. Using again that $L(x)\leq G(x)\leq U(x)$ we see that $|\bar{G}_{A,B}(x)-G(x)|\leq U(x)-B < G(x_b)-B+\varepsilon \leq |G(x_b)-h(x_b)|+\varepsilon$
for small enough $\varepsilon$, showing that $|\bar{G}_{A,B}(x)-G(x)|\leq \|h-G\|$.

Similarly, if $\bar{G}(x)<A$ we conclude that  $|\bar{G}_{A,B}(x)-G(x)|\leq \|h-G\|$. 

Finally, assume that
$U(x)>L(x)$ and $\bar{G}(x)\in [A,B]$.
Since $h$ is nonincreasing we have that $h(x_a)\geq h(x_b)$ and, consequently, $$\|h-G\|\geq \max(|h(x_a)-G(x_a)|,|h(x_b)-G(x_b)|)\geq \frac{G(x_b)-G(x_a)}2\geq |\bar{G}_{A,B}(x)-G(x)|-2\varepsilon$$
for $\varepsilon$ small enough. This completes the proof of (\ref{bestappr}).
\vspace{3mm}

To check continuity of $U$
note that for $0\leq y<x\leq 1$ $U(y)=\max(U(x),\max_{y\leq z\leq x} G(z))$. Now, given $\varepsilon>0$ we can fix $\delta>0$ such that
$|G(x)-G(y)|\leq \varepsilon$ whenever $|y-x|\leq \delta$. But then $|U(y)-U(x)|\leq \varepsilon$ if $|y-x|\leq \delta$, proving continuity
of $U$. $L$ can be handled similarly. As a consequence we see that $\bar{G}$ and $\bar{G}_{A,B}$ are also continuous. 
\vspace{3mm}

Now, to prove the first equality in the statement
we take $x\in [0,1]$ and consider first the case $x\in \mathcal{T}_1$. Note that, necessarily, $U(x)\geq B$, $U(x)-B\geq B-L(x)$ and $\bar{G}_{A,B}(x)=B$.

\noindent
$\bullet$\hspace{3mm}If $G(x)\geq B$ then $|G(t)-\bar{G}_{A,B}(x)|=G(x)-B$. 

\noindent
$\bullet$\hspace{3mm}Assume, on the contrary, that $G(x)<B$. Set $x_+=\inf\{ y\leq x:\, G(y)=U(x)\}$. By continuity,
$G(x_+)=U(x)=U(x_+)$. 

Now, if $\bar{G}(x_+)\geq B$ then $G(x_+)-B=U(x)-B\geq B-L(x)\geq B-G(x)=|G(x)-\bar{G}_{A,B}(x)|$. If, on the contrary, 
$\bar{G}(x_+)<B$, then there exists $x'\in[x,x_+]$ such that $\bar{G}(x')\in(A,B)$. But we must have $U(x')=U(x)=U(x_+)$  and $L(x')<L(x)$ and,
consequently, we have that $$|G(x)-\bar{G}_{A,B}(x)|=B-G(x)\leq B-L(x)\leq \frac {U(x)-L(x)}2 < \frac {U(x')-L(x')}2.$$ 

Summarizing, we see that
\begin{equation}\label{bd1}
\max_{\bar{G}(x)\geq B}|G(x)-\bar{G}_{A,B}(x)|\leq \max\left(\max_{ \bar{G}(x)\geq B} (G(x)-B),{\textstyle \frac 1 2}\max_{A\leq \bar{G}(xt)\leq B} (U(x)-L(x))\right).
\end{equation}

Similarly, 
\begin{equation}\label{bd2}
\max_{\bar{G}(x)\leq A}|G(x)-\bar{G}_{A,B}(x)|\leq \max\left(\max_{ \bar{G}(x)\leq A} (A-G(x)),{\textstyle \frac 1 2}\max_{A\leq \bar{G}(x)\leq B} (U(x)-L(x))\right)
\end{equation}
and, obviously, if $\bar{G}(x)\in[A,B]$ then $\bar{G}_{A,B}(x)=\bar{G}(x)$ and $|G(x)-\bar{G}_{A,B}(x)|\leq \frac 1 2 (U(x)-L(x))$, which implies that
\begin{equation}\label{bd3}
\max_{A\leq \bar{G}(x)\leq B}|G(x)-\bar{G}_{A,B}(t)|\leq {\textstyle \frac 1 2}\max_{A\leq \bar{G}(x)\leq B} (U(x)-L(x)).
\end{equation}
Now combining (\ref{bd1}), (\ref{bd2}) and (\ref{bd3}) we see that 
$$\|G-\bar{G}_{A,B}\|\leq \max\left(\max_{ \bar{G}(x)\geq B} (G(x)-B),\max_{\bar{G}(x)\leq A} (A-G(x)),{\textstyle \frac 1 2}\max_{A\leq \bar{G}(x)\leq B} (U(x)-L(x))\right).$$
\vspace{3mm}

Assume now that $x_0$ is such that $\bar{G}(x_0)\geq B$. Then $\bar{G}_{A,B}(x_0)=B$ and $G(x_0)-B\leq |G(x_0)-\bar{G}_{A,B}(x_0)|$.  This implies 
$\max_{ \bar{G}(x)\geq B} (G(x)-B)\leq \|G-\bar{G}_{A,B}\|$. 

Similarly, $\max_{ \bar{G}(t)\leq A} (A-G(x))\leq \|G-\bar{G}_{A,B}\|$.
\vspace{3mm}

Finally, suppose $x_0$ is such that $\bar{G}(x_0)\in [A,B]$ and $$U(x_0)-L(x_0)=\max_{\bar{G}(x)\in[1,B]} (U(x)-L(x))\geq \max\left(\max_{ \bar{G}(x)\geq B} (G(x)-B), 
\max_{ \bar{G}(x)\leq A} (A-G(x))\right).$$ 

\noindent
$\bullet$\hspace{3mm}If $U(x_0)=L(x_0)$ then 
$$\|G-\bar{G}_{A,B}\|= \max\left(\max_{ \bar{G}(x)\geq B} (G(x)-B),\max_{\bar{G}(x)\leq A} (A-G(x)),{\textstyle \frac 1 2}\max_{A\leq \bar{G}(x)\leq B} (U(x)-L(x))\right)=0.$$

\noindent
$\bullet$\hspace{3mm}If $U(x_0)>L(x_0)$ then we set $x_+=\inf\{y\in[x_0,1]:\, G(y)=U(x_0)  \}$. Then $U(y)=U(x_0)$ for $y\in [x_0,x_+]$ and
$$G(x_+)=U(x_+)=U(x_0).$$ Set $x_+=\sup\{y\in [0,x_0]:\, G(y)=L(x_0)\}$. We have $L(y)=L(x_0)=G(x_-)$ for $y\in[x_-,x_0]$. 
We claim that 
\begin{equation}\label{claim1}
L(y)=L(x_0) \quad  \mbox{ for }y\in [x_0,x_+].
\end{equation}
To check (\ref{claim1}) note that, if $\bar{G}(x_0)>A$ and (\ref{claim1}) fails then we could find $y\in [x_0,x_+]$ with $L(y)<L(x_0)$, $\bar{G}(y)\in (A,B]$ and $U(y)-L(y)>U(x_0)-L(x_0)$,
while if $\bar{G}(x_0)=A$ and (\ref{claim1}) fails then $G(y)<L(x_0)$ for some $y\in (x_0,x_+)$, $\bar{G}(y)<A$ and $A-L(y)>A-L(x_0)=\frac 1 2(U(x_0-L(x_0)$, against the assumption on $x_0$.

Hence, from (\ref{claim1}) we conclude that $\bar{G}(x_+)=\bar{G}(x_0)\in[A,B]$ and $|G(x_+)-\bar{G}_{A,B}(x_+)|= \frac 1 2 (U(x_0)-L(x_0))$, showing that $\frac 1 2 (U(x_0)-L(x_0))\leq 
\|G-\bar{G}_{A,B}\|$. Combining the last estimates we see that the first equality in (\ref{rep1}) holds.
\vspace{3mm}

For the second identity we note that arguing as above we see that $U(x_0)-L(x_0)=G(x)-G(y)$ for some $(y,x)\in \mathcal{T}_3$ if $\bar{G}(x_0)\in [A,B]$. Assume, on the other hand, that $(y_0,x_0)\in \mathcal{T}_3$
satisfies $$\frac  1 2 (G(x_0)-G(y_0))\geq  \max\left(\max_{ \bar{G}(x)\geq B} (G(x)-B),\max_{\bar{G}(x)\leq A} (A-G(x))\right).$$ 

\noindent
$\bullet$\hspace{3mm}We consider first the case $\frac 1 2 (G(y_0)+G(x_0))\in (A,B)$.

We claim that $U(x_0)=G(x_0)$ since, otherwise, there exists $x'>x_0$ such that $\frac 1 2 (G(y_0)+G(x'))\in (A,B)$ and $G(x')>G(x_0)$ and this would imply $G(x')-G(y_0)>G(x_0)-G(y_0)$, against the assumption.

Similarly, we see that $G(y_0)=L(x_0)$. 

Furthermore, $L(x)=L(y_0)$ for $x\in [y_0,x_0]$. If $G(x_0)<U(x_0)$ then there exists $x'>x_0$ such that $\frac 1 2 (G(y_0)+G(x'))\in (A,B)$ and $G(x')>G(x_0)$,
but then $G(x')-G(y_0)>G(x_0)-G(y_0)$, contradicting maximality of $(y_0,x_0)$. Similarly we see that $G(y_0)=L(y_0)$ and also that $L(x)=L(y_0)$ for $x\in [y_0,x_0]$. Hence,
$G(x_0)-G(y_0)=U(x_0)-L(x_0)$ and $\bar{G}(x_0)\in(A,B)$. 

\noindent
$\bullet$\hspace{3mm}In the case $\frac 1 2 (G(y_0)+G(x_0))=B$ we have that necessarily $G(x_0)\geq B$ and, arguing as above, we see that $G(y_0)=L(y)$ for all
$y\in [y_0,x_0]$. This implies that $\bar{G}(x_0)\geq B$ and $\frac 1 2 (G(x_0)-G(y_0))=G(x_0)-B$. 

\noindent
$\bullet$\hspace{3mm}Arguing similarly for the case  $\frac 1 2 (G(y_0)+G(x_0))=A$ we conclude that the second
equality in (\ref{rep1}) holds.
\hfill $\Box$

\begin{Remark} \label{notaultima}
The sets of optimizers within $\mathcal{T}_1, \mathcal{T}_2$ and $\mathcal{T}_3$ in Lemma \ref{representation1} play an important role in the next results. 
For convenience, we denote ${T}_1=\{x_0\in \mathcal{T}_1:\, G(x_0)-B=\|G-\bar{G}_{A,B}\|\}$, ${T}_2=\{x_0\in \mathcal{T}_2:\, A-G(x_0)=\|G-\bar{G}_{A,B}\|\}$
and ${T}_3=\{(y_0,x_0)\in \mathcal{T}_3:\, \frac 1 2 (G(x_0)-G(y_0))=\|G-\bar{G}_{A,B}\|\}$. A look at the proof of Lemma \ref{representation1} shows that
if $x_0\in T_1$ then $G$ has a local maximum at $x_0$ and a local minimum if $x_0\in T_2$. Also, if $(y_0,x_0)\in T_3$ then $G$ has a local maximum at $x_0$ and a local minimum 
at $y_0$.
\end{Remark}

\medskip
Our next result addresses the directional differentiability of the functional $G\to \|G-\bar G_{A,B}\|$ that appeared in the last theorem. This kind of result typically allows to obtain efficiency and asymptotic distributional behaviour of  functionals  in the statistical setting {(see e.g. \cite{Carcamo})}. In fact it allows to prove the Central Limit Theorem for the statistical functional $d_K(F_0,R_\alpha(F_n))$ (see Theorem 4.1 in \cite{Hristo}).

\medskip
\begin{Theorem}\label{asymptotic1}
Assume $G,J:[0,1]\to \mathbb{R}$ are continuous functions and $r_n>0$ is a sequence of real numbers such that $r_n\to \infty$. 
Define $G_n=G+\frac J{r_n}$ and consider $\bar{G}, \bar{G}_{A,B}$ as in Theorem \ref{representation1} and $\bar{G}_{A,B,n}$ built in the 
same way as $G_{A,B}$ but from $G_n$. Assume further that $T_1, T_2$ and $T_3$ are as in Remark \ref{notaultima} and that
there is no $x\in {T}_1$ with $\bar{G}(x)=B$, no $x\in {T}_2$ with $\bar{G}(x)=A$ and no $(y,x)\in {T}_3$ with $\frac 1 2 (G(x)+G(y))\in \{A,B\}$.
Then
$$r_n(\|G_n-\bar{G}_{A,B,n}\|-\|G-\bar{G}_{A,B}\|)\to \max\left(\max_{x\in {T}_1} J(x), \max_{t\in {T}_2} (-J(x)),\frac 1 2 \max_{(y,x)\in {T}_3} (J(x)-J(y)) \right).$$
\end{Theorem}

\medskip
\noindent \textbf{Proof.} We use the notation $U,L$ from Theorem \ref{representation1} and write $U_n,L_n,\bar{G}_n, T_{n,i}$ for the corresponding objects coming from $G_n$.
Observe that $\|U_n-U\|\leq \|J\|/r_n\to 0$ and, similarly, $\|\bar{G}_n-\bar{G}\|\to 0$. Assume that $x\in {T}_1$. By assumption and the last convergence we have that 
$\bar{G}_n(x)>B$ for large enough $n$ and, therefore,
$\|G_n-\bar{G}_{A,B,n}\|\geq (G_n(t)-B)$. But this implies
$$r_n(\|G_n-\bar{G}_{A,B,n}\|-\|G-\bar{G}_{A,B}\|)\geq r_n((G_n(x)-B)-(G(x)-B))=J(x).$$
Arguing similarly for ${T}_2$ and ${T}_3$ we conclude that
\begin{eqnarray}\label{cotinf}
\lefteqn{\liminf r_n(\|G_n-\bar{G}_{A,B,n}\|-\|G-\bar{G}_{A,B}\|)}\hspace*{3cm}\\
\nonumber
&\geq &\max\left(\max_{x\in {T}_1} J(x), \max_{x\in {T}_2} (-J(x)),\frac 1 2 \max_{(y,x)\in {T}_3} (J(x)-J(y)) \right).
\end{eqnarray}
For the upper bound assume $x_n\in {T}_{n,1}$ (that is, $x_n\in \mathcal{T}_{n,1}$ such that $G_n(x_n)-B=\|G_n-\bar{G}_{A,B,n}\|)$. By compactness,
taking subsequences if necessary, we can assume that $x_n\to x_0$ for some $x_0\in [0,1]$ with $\bar{G}(x_0)\geq B$ and 
$G(x_0)-B=\|G-\bar{G}_{A,B}\|$. But this means that $x_0\in {T}_1$. Hence, by assumption $G(x_0)>B$ and, consequently,
$G(x_n)>B$ for large enough $n$. In this case $\|G-\bar{G}_{A,B}\|\geq (G(x_n)-B)$, which implies that
$$r_n(\|G_n-\bar{G}_{A,B,n}\|-\|G-\bar{G}_{A,B}\|)\leq r_n((G_n(x_n)-B)-(G(x_n)-B))=J(x_n)\to J(x_0).$$
With the same argument applied to ${T}_2$ and ${T}_3$ we conclude that 
\begin{eqnarray}\label{cotsup}
\lefteqn{\limsup r_n(\|G_n-\bar{G}_{A,B,n}\|-\|G-\bar{G}_{A,B}\|)}\hspace*{3cm}\\
\nonumber
&\leq &\max\left(\max_{x\in {T}_1} J(x), \max_{x\in {T}_2} (-J(x)),\frac 1 2 \max_{(y,x)\in {T}_3} (J(x)-J(y)) \right)
\end{eqnarray}
and complete the proof.
\hfill $\Box$
\vspace{3mm}

Specializing the last results for  $G(x)=f(x)-Lx$, where $f$ is nondecreasing, $L\geq 1$  a  constant, and $A=1-L, B=0$, we can obtain a first result on the directional differentiability of the functional $f\to \|f-\tilde f_L\|$ considered in Section \ref{Approximation}. Note that now, recovering the notation in that section, the relevant sets are $\mathcal{T}_1, \mathcal{T}_2$ and $\mathcal{T}_3$ as defined in Theorem \ref{fcont}, and ${T}_1=\{x_0\in \mathcal{T}_1:\, f(x_0)-Lx_0=\|f-\tilde f_L\|\}$, ${T}_2=\{x_0\in \mathcal{T}_2:\, 1-L+Lx_0-f(x_0)=\|f-\tilde f_L\|\}$
and ${T}_3=\{(y_0,x_0)\in \mathcal{T}_3:\, \frac 1 2 (f(x_0)-f(y_0)-L(x_0-y_0))=\|f-\tilde f_L\|\}$. Theorem \ref{asymptotic1} translates then to the following immediate corollary.

\begin{Corollary}[Directional differentiability.]\label{derivability}
Let $f, f_n:[0,1]\to \mathbb{R}$ be  nondecreasing  functions, $r_n>0$  a sequence of real numbers such that $r_n\to \infty$ and $r_n(f_n-f) \to J$ pointwise, where $J:[0,1]\to \mathbb{R}$ is a continuous function. Assume further that $f$ is continuous, that $T_1, T_2$ and $T_3$ are as above and that
there is no $x\in {T}_1$ with $\gamma_L(x)=0$, no $x\in {T}_2$ with $\gamma_L(x)=1-L$ and no $(y,x)\in {T}_3$ with $\frac 1 2 (f(x)+f(y)-L(x+y))\in \{1-L,0\}$. Let  $\tilde f_{n,L}, \tilde f_L$ respectively denote the best $L_\infty$-approximations to $f_n$ and $f$ by Lipschitz-continuous functions $h:[0,1]\to \mathbb R$ with $\|h\|_\lip\leq L$ and verifying $h(0)=0, h(1)=1$, as in Theorem \ref{final}. Then
$$r_n(\|f_n-\tilde f_{L,n}\|-\|f-\tilde f_L\|)\to \max\left(\max_{x\in {T}_1} J(x), \max_{x\in {T}_2} (-J(x)),\frac 1 2 \max_{(y,x)\in {T}_3} (J(x)-J(y)) \right).$$
\end{Corollary}

\end{document}